\documentclass[12pt]{article}

\usepackage{acronym}
\usepackage{amsfonts}
\usepackage{amsmath,amssymb,bm}
\usepackage{graphicx}
\usepackage{subfigure}
\usepackage{multirow}

\newtheorem{theorem}{Theorem} %[section]

\newtheorem{definition}{Definition}

\def\be{\begin{Example}}
\def\ee{\end{Example}}
\def\bt{\begin{theorem}}
\def\et{\end{theorem}\bigskip}
\def\bl{\begin{Lemma}}
\def\el{\end{Lemma}\bigskip}
\def\ep{\end{Proposition}\bigskip}
\def\bp{\begin{Proposition}}
\def\bd{\begin{definition}}
\def\ed{\end{definition}}

\begin{document}
\title{\bf An adaptive gradient method for computing generalized tensor eigenpairs
% \thanks{This work was supported in part by the National Natural Science Foundation of China (No.61262026), NCET Programm of the Ministry of Education (NCET 13-0738), JGZX programm of Jiangxi Province (20112BCB23027), Natural Science Foundation of Jiangxi Province (20132BAB201026), science and technology programm of Jiangxi Education Committee (LDJH12088). }
 }
 \author{
 %Liqun Qi\footnote{Department of Applied
% Mathematics, The Hong Kong Polytechnic University, Hung Hom, Kowloon, Hong Kong.
% E-mail: maqilq@polyu.edu.hk.}\quad
 Gaohang Yu \thanks{\small
  School of Mathematics and Computer Sciences,
Gannan Normal University, Ganzhou, 341000, China. E-mail:
maghyu@163.com}
\quad Zefeng Yu \thanks{\small
  School of Mathematics and Computer Sciences,
Gannan Normal University, Ganzhou, 341000, China. E-mail:
yzf\_2000@sina.com}
\quad Yi Xu \thanks{\small
 Department of Mathematics, Southeast University, China. E-mail: yi.xu1983@gmail.com}
 \quad Yisheng Song \thanks{\small
 School of Mathematics and Information Science and Henan Engineering Laboratory for Big Data Statistical Analysis and Optimal Control, Henan Normal University, China. E-mail: songyisheng1@gmail.com}
  }

\date{\today}
\maketitle

%{\large

\begin{abstract}
High order tensor arises more and more often in signal processing,
data analysis, higher-order statistics, as well as imaging sciences. In this paper, an adaptive gradient (AG) method is presented for generalized tensor eigenpairs. Global convergence and linear convergence rate are established under some suitable conditions. Numerical results are reported to illustrate the efficiency of the proposed method. Comparing with the GEAP method, an adaptive shifted power method proposed by Tamara G. Kolda and Jackson R. Mayo [SIAM J. Matrix Anal. Appl., 35 (2014), pp. 1563-1581],  the AG method is much faster and could reach the largest eigenpair with a higher probability.

{\bf Keywords:} Higher order Tensor, Eigenvalue, Eigenvector, Gradient method, Power method.

\end{abstract}

\section{Introduction}

A $m$th-order $n$-dimensional real tensor $\mathcal{A}$ consists of
$n^m$ entries in real numbers:
\[
\mathcal{A}=(a_{i_1i_2\cdots
i_m}),\,\, a_{i_1i_2\cdots i_m}\in \mathbb{R},\,\, 1\leq
i_1,i_2,\ldots,i_m\leq n.
\]
$\mathcal{A}$  is called \textit{symmetric} if the value of
$a_{i_1i_2\cdots i_m}$ is invariant under any
permutation of its indices $i_1,i_2,\ldots,i_m$.
Recall the definition of tensor product, $\mathcal{A}x^{m-1}$ is a vector in $\mathbb{R}^n$ with its $i$th
component as
\begin{equation}\label{e-produdefn}
(\mathcal{A}x^{m-1})_i=\sum_{i_2,\ldots,i_m=1}^n a_{i i_2 \cdots i_m}x_{i_2}\cdots
x_{i_m}.
\end{equation}
 A real symmetric tensor $\mathcal{A}$ of order $m$ dimension $n$ uniquely defines a
$m$th degree homogeneous polynomial function $h$ with  real
coefficient by
\begin{align}\label{f_poly}
h(x):=\mathcal{A}x^m=x^T(\mathcal{A}x^{m-1})=\sum_{i_1,\ldots, i_m=1}^n a_{i_1 \cdots i_m}x_{i_1}\cdots
x_{i_m}.
\end{align}
We call that the tensor $\mathcal{A}$ is positive definite if $\mathcal{A}x^m>0$ for all $x\neq 0$.

In 2005, Qi \cite{Qi05} and Lim \cite{Lim05} proposed the definition of eigenvalues and eigenvectors for higher order tenors, independently.
Furthermore, in \cite{CPZ09}, these definitions were unified by Chang, Person and Zhang. Let $\mathcal{A}$ and $\mathcal{B}$ be real-valued, $m$th-order $n$-dimensional symmetric tensors. Assume further that $m$ is even and $\mathcal{B}$ is positive definite. we call $(\lambda, x) \in \mathbb{R} \times \mathbb{R}^n \backslash \{0\}$ is a {\bf generalized eigenpair} of $(\mathcal{A},\mathcal{B})$ if
\begin{equation}\label{generalized eigenpair}
\mathcal{A}x^{m-1}=\lambda\mathcal{B}x^{m-1}.
\end{equation}

    When the tensor $\mathcal{B}$ is an identity tensor $\mathcal{\varepsilon}$ such that $\mathcal{\varepsilon}x^{m-1}=\|x\|^{m-2}x$ for all $x \in \mathbb{R}^n$ \cite{CPZ09}, the eigenpair reduces to $Z$-eigenpair \cite{Qi05,Lim05} which is defined as a pair $(\lambda, x) \in \mathbb{R} \times \mathbb{R}^n \backslash \{0\}$ satisfying
  \begin{equation}\label{Zeigenpair}
  \mathcal{A}x^{m-1}=\lambda x  \ \ \mbox{and} \ \ \|x\|^2=1.
  \end{equation}
 Another special case is that when $\mathcal{B}=\mathcal{I}$ with $(\mathcal{I})_{i_1 i_2\cdots i_m}=\delta_{i_1 i_2\cdots i_m}$ \cite{CPZ09}, the real scalar $\lambda$ is called an $H$-eigenvalue and the real vector $x$ is the associated $H$-eigenvector of the tensor $\mathcal{A}$ \cite{Qi05}.

In the last decade, tensor eigenproblem has received much attention in the literature \cite{CQZ13,DW15,LQY13-LAA,SQ15,SY15,YYZ14,YY11,YY10}, which has numerous applications in magnetic resonance imaging \cite{SS08,BV08,QWW08,QHW09,QYW,QYX,CDHS13}, image analysis \cite{ZZP13}, data fitting \cite{PLV05,PLV09}, quantum information \cite{NQB14}, automatic control \cite{NQW08}, higher order Markov chains \cite{LN13,CZ13}, spectral graph theory \cite{LQY13,XC13}, multi-label learning \cite{SJY08}, and so on. In \cite{QYW}, a positive semidefinite diffusion tensor (PSDT) model was proposed to approximate the apparent diffusion coefficient (ADC) profile for high-order diffusion tensor imaging, where the smallest Z-eigenvalue need to be nonnegative to guarantee the positive definiteness of the diffusivity function. Based on all of the Z-eigenvalues, a generalized fractional anisotropy (GFA) was proposed to characterize the anisotropic diffusion profile for PSDT. GFA is rotationally invariant and independent from the choice of the laboratory coordinate system. In automatic control \cite{NQW08}, the smallest eigenvalue of tensors could reflect the stability of a nonlinear autonomous system. In \cite{QYX}, the principal Z-eigenvectors can depict the orientations of nerve fibers in the voxel of white matter of human brain. Recently, a higher order tensor vessel tractography was proposed for segmentation of vascular structures, in which the principal directions of a 4-dimensional tensor were used in vessel tractography approach \cite{CetinUnal15}.

In general, it is NP-hard to compute eigenvalues of a tensor \cite{HillarLim2013}. In \cite{QWW}, a direct method to calculate all of Z-eigenvalues was proposed for two and three dimensional symmetric tensors. For general symmetric tensors, a shifted higher order power method was proposed for computing Z-eigenpairs in \cite{KoldaMayo11}. Recently, in \cite{KoldaMayo14}, an adaptive version of higher order power method was presented for generalized eigenpairs of symmetric tensor. In order to guarantee the convergence of power method, they need a shift to force the objective to be (locally) concave/convex. In this case, the power method is a monotone gradient method with unit-stepsize. By using fixed-point analysis, linear convergence rate is established for the shifted higher order power method \cite{KoldaMayo11}. However, similarly to the case of Matrix, when the largest eigenvalue is close to the second dominant eigenvalue, the convergence of power method will be very slow \cite{CW14}.

In the recent years, there are various optimization approaches were proposed for tensor eigenvalue problem \cite{Han2013,HLQ13,HCD15,HCD15-PJO, NQ15-GlobalOpt,ZN15-PJO}.  In \cite{Han2013}, Han proposed an unconstrained optimization model for computing generalized eigenpair of symmetric tensors. By using BFGS method to solve the unconstrained optimization, the sequence will be convergent superlinearly. A subspace projection method was proposed in \cite{HCD15} for Z-eigenvalues of symmetric tensors. Recently, in \cite{HCD15-PJO}, Hao, Cui and Dai proposed a trust region method for Z-eigenvalues of symmetric tensor and the sequence enjoys a locally quadratic convergence rate.  In \cite{ZN15-PJO}, Ni and Qi employed Newton method for the KKT system of optimization problem, and obtained a quadratically convergent algorithm for finding the largest eigenvalue of a nonnegative homogeneous polynomial map. In \cite{ChenQi15}, an inexact steepest descent method was proposed for computing eigenvalues of large scale Hankel tensors. Since nonlinear optimization methods may stop at a local optimum, a sequential semi-definite programming method was proposed by Hu et al. \cite{HHQ13} for finding the extremal Z-eigenvalues of tensors. Moreover, in \cite{CDN14}, a Jacobian semi-definite relaxation approach was presented to compute all of the real eigenvalues of symmetric tenors.

In practice, one just need to compute extremal eigenvalues or all of its local maximal eigenvalues, for example in MRI \cite{QYW, QYX}.  On the other hand, when the order or the dimension of a tensor grows larger, the optimization problem will become large-scale or huge-scale. Therefore, we would like to investigate one simple and low-complexity method for finding tensor eigenpairs.
In this paper, we consider an adaptive gradient method for solving the following nonlinear programming problem:
\begin{equation}\label{max-optimization-problem}
\max f(x)=\frac{\mathcal{A}x^m}{\mathcal{B}x^m} \ \ \mbox{subject to} \ \ x \in \mathbb{S}^{n-1}.
\end{equation}
 Where $\mathbb{S}^{n-1}$ denote the unit sphere, i.e., $\mathbb{S}^{n-1}=\{x\in \mathbb{R}^n | \|x\|^2=1\}$, $\|\cdot\|$ denotes the Euclidean norm.
By some simple calculations, we can get its gradient and Hessian, as follows \cite{KoldaMayo14,ChenQi15}:
\begin{equation}\label{gradient}
g(x)\equiv \nabla f(x)=\frac{m}{\mathcal{B}x^m}(\mathcal{A}x^{m-1}-\frac{\mathcal{A}x^m}{\mathcal{B}x^m}\mathcal{B}x^{m-1}),
\end{equation}
and its Hessian is
\begin{align*}\label{Hessian}
H(x)\equiv & \nabla^2 f(x)\\
&=\frac{m(m-1)\mathcal{A}x^{m-2}}{\mathcal{B}x^m}-\frac{m(m-1)\mathcal{A}x^m\mathcal{B}x^{m-2}+m^2(\mathcal{A}x^{m-1}\circledcirc\mathcal{B}x^{m-1})}{(\mathcal{B}x^m)^2}\\
&+\frac{m^2\mathcal{A}x^m(\mathcal{B}x^{m-1}\circledcirc\mathcal{B}x^{m-1})}{(\mathcal{B}x^m)^3},
\end{align*}
where $x\circledcirc y= xy^T+yx^T$, and $\mathcal{A}x^{m-2}$ is a matrix with its component as $$(\mathcal{A}x^{m-2})_{i_1i_2}=\sum_{i_3,\ldots,i_m=1}^n a_{i_1 i_2i_3 \cdots i_m}x_{i_3}\cdots
x_{i_m}\ \ \mbox{for all}\ \ i_1,i_2=1,\cdots,n.$$

According to (\ref{gradient}), we can derive an important property for the nonlinear programming problem (\ref{max-optimization-problem}) that the gradient $g(x)$ is located in the tangent plane of  $\mathbb{S}^{n-1}$ at $x$ \cite{ChenQi15}, since
 \begin{equation}\label{xtgx}
x^Tg(x)=\frac{m}{\mathcal{B}x^m}(x^T\mathcal{A}x^{m-1}-\frac{\mathcal{A}x^m}{\mathcal{B}x^m}x^T\mathcal{B}x^{m-1})=0.
\end{equation}

Let $\overline{x}$ is a constrained stationary point of (\ref{max-optimization-problem}), i.e., that
$\langle g(\overline{x}),x-\overline{x}\rangle \le 0 \ \ \mbox{for all}\  x \in \mathbb{S}^{n-1}.$
Then we can claim that {\bf every constrained stationary point} of (\ref{max-optimization-problem}) {\bf must be a stationary point} of $f(x)$ since
$\langle g(\overline{x}),x\rangle \le 0$ should be hold for all $ x \in \mathbb{S}^{n-1}.$  Otherwise, if $\|g(\overline{x})\| \neq 0$, we could choose $x=\frac{g(\overline{x})}{\|g(\overline{x})\|}$, and then $\|g(\overline{x})\| \le 0$.

Suppose $x \in \mathbb{S}^{n-1}$ and denote $\lambda = \frac{\mathcal{A}x^m}{\mathcal{B}x^m}$. By $g(x)=0$, we know that any KKT point of (\ref{max-optimization-problem}) will be a solution of the system of equations (\ref{generalized eigenpair}).
 Before end of this section, we would like to state the following theorem and its proof is omitted.

 \begin{theorem}   If the gradient $g(x)$ at $x$ vanishes, then $\lambda=f(x)$ is a generalized eigenvalue and the vector $x$ is the associated generalized eigenvector.
\end{theorem}

The rest of this paper is organized as follows. In the next section, we introduce some existed gradient methods for tensor eigenvalue problems.
In Section 3, based on a curvilinear search scheme, we present a inexact gradient method. Then, we establish its global convergence and linear convergence results under some suitable assumptions. Section 4 provides numerical experiments to show the efficiency of our gradient method. Finally, we have a conclusion section.

\section{Some exist gradient methods for tensor eigenpairs}
\subsection{Gradient method with fixed stepsize--power method}
 The symmetric higher-order power method (S-HOPM) was introduced by De Lathauwer, De Moor, and Vandewalle \cite{LMV00} for solving the following optimization problem:
\begin{equation}\label{MaxAxm}
\max f(x)=\mathcal{A}x^m \ \ \mbox{subject to} \ \ x \in \mathbb{S}^{n-1}.
\end{equation}
%or
%\begin{equation}\label{MinAxm}
%\min \mathcal{A}x^m \ \ \mbox{subject to} \ \ x \in \mathbb{S}^{n-1}.
%\end{equation}
This problem is equivalent to finding the largest Z-eigenvalue of $\mathcal{A}$ \cite{Qi05} and is related to finding the best symmetric rank-1 approximation of a symmetric tensor $\mathcal{A}\in \mathbb{S}^{[m,n]}$ \cite{LMV00}.

\noindent\rule{\textwidth}{1pt}
\underline{\textbf{Algorithm 1: Symmetric higher order power Method (S-HOPM) \cite{LMV00} }} \\
Given a tensor $\mathcal{A} \in \mathbb{S}^{[m,n]}$,
an initial unit iterate $x_0$. Let $\lambda_0=\mathcal{A}x_0^m.$ \\
\textbf{for} $k=0,1,\ldots$ \textbf{do}\\
%calculate a search direction $d(x_k)$.\\
\textbf{1:} $\hat{x}_{k+1}\leftarrow \mathcal{A}x_k^{m-1}$\\
\textbf{2:} $x_{k+1}\leftarrow \hat{x}_{k+1}/\|\hat{x}_{k+1}\|$\\
\textbf{3:} $\lambda_{k+1}=\mathcal{A}x_{k+1}^m$ \\
\textbf{End for}\\
\noindent\rule{\textwidth}{1pt}

 The cost per iteration of power method is $O(mn^m)$, mainly for computing $\mathcal{A}x_k^{m-1}$. Let $g(x_k)=\nabla f(x_k)=\frac{1}{m}\mathcal{A}x_k^{m-1}$. Set $d_k=g(x_k)-x_k$, $\hat{x}_{k+1}=x_k+d_k$, then the main iteration could be reformulated as $x_{k+1}=\frac{x_k+d_k}{\|x_k+d_k\|}$, which is a projected gradient method with unit-stepsize. Kofidis and Regalia \cite{KR02} pointed out that S-HOPM method can not guarantee to converge. By using convexity theory, they show that S-HOPM method could be convergent for even-order tensors under the convexity assumption on $\mathcal{A}x^m$. For general symmetric tensors, a shifted S-HOPM (SS-HOPM) method was proposed by Kolda and Mayo \cite{KoldaMayo11} for computing Z-eigenpairs. One shortcoming of SS-HOPM is that its performance depended on choosing an appropriate shift. Recently, Kolda and Mayo extended SS-HOPM for computing generalized tensor eigenpairs, called GEAP method which is an adaptive, monotonically convergent, shifted power method for generalized tensor eigenpairs (\ref{generalized eigenpair}). They showed that GEAP method is much faster than the SS-HOPM method due to its adaptive shift choice.

\noindent\rule{\textwidth}{1pt}
\underline{\textbf{Algorithm 2: GEAP Method \cite{KoldaMayo14} }} \\
Given tensors $\mathcal{A} \in \mathbb{S}^{[m,n]}$ and $\mathcal{B}\in \mathbb{S}_+^{[m,n]}$, and an initial guess $\hat{x}_0$. Let $\beta=1$ if we want to find the local maxima; otherwise, let $\beta=-1$ for seeking local minima. Let $\tau>0$ be the tolerance on being positive/negative definite.\\
\textbf{for} $k=0,1,\ldots$ \textbf{do}\\
\textbf{1:} Precompute $\mathcal{A}x_k^{m-2}$,$\mathcal{B}x_k^{m-2}$,$\mathcal{A}x_k^{m-1}$,$\mathcal{B}x_k^{m-1}$,$\mathcal{A}x_k^{m}$, $\mathcal{B}x_k^{m}$ \\
\textbf{2:} $\lambda_{k}=\mathcal{A}x_{k}^m/\mathcal{B}x_{k}^m$ \\
\textbf{3:} $H_k\leftarrow H(x_k)=\nabla^2 f(x_k)$\\
\textbf{4:} $\alpha_k\leftarrow \beta \max\{0,(\tau-\lambda_{min}(\beta H_k))/m$\\
\textbf{5:} $\hat{x}_{k+1}\leftarrow \beta(\mathcal{A}x_k^{m-1}-\lambda_k \mathcal{B}x_k^{m-1}+(\alpha_k+\lambda_k)\mathcal{B}x_k^m x_k) $\\
\textbf{6:} $x_{k+1}\leftarrow \hat{x}_{k+1}/\|\hat{x}_{k+1}\|$\\
\textbf{End for}\\
\noindent\rule{\textwidth}{1pt}

%Noticed that the objective function used in GEAP method is $f(x)=\frac{\mathcal{A}x^m}{\mathcal{B}x^m}\|x\|^m+\alpha \|x\|^m$. And the main iteration is also could be reformulated as  $x_{k+1}=P_{\mathbb{S}^{n-1}}(x_k+d_k)$ with $d_k=g_k-x_k$. Here $P_{\mathbb{S}^{n-1}}(x)=\frac{x}{\|x\|}$ is a project operate. The cost per iteration of GEAP method is $O(n^m/m!)$ for symmetry tensor.

In \cite{NQZ10}, Ng, Qi and Zhou proposed a power method for finding the largest H-eigenvalue of irreducible nonnegative tensors. It is proved in \cite{CPZ11} that NQZ's power method is convergent for primitive nonnegative tensors. Further, Zhang et. al \cite{ZQ12,ZQX12} established its linear convergence result and presented some updated version for essentially positive tensors and weakly positive tensors, respectively.
 %Liu et. al \cite{LZI10} proposed a shifted power method for any irreducible nonnegative tensor.
 However, similarly to the case of Matrix, when the largest eigenvalue is close to the second dominant eigenvalue, the convergence of power method will be very slow \cite{CW14}.

\subsection{Gradient method with optimal stepsize}

In \cite{HCD15}, Hao, Cui and Dai proposed a sequential subspace projection method (SSPM) for Z-eigenvalue of symmetric tensors. In each iteration of SSPM method,  one need to solve the following 2-dimensional subproblem:

\begin{equation}\label{2-dimenional subproblem}
\max_{x \in span\{x_k, \mathcal{A}x_k^{m-1}\}} f(x)=\mathcal{A}x^m \ \ \mbox{subject to} \ \ x \in \mathbb{S}^{n-1}.
\end{equation}
Let $g_k\triangleq \mathcal{A}x_k^{m-1}$. The point in $\mathbb{S}^{n-1}\cap span\{x_k, g_k\}$ can be expressed as
$$x(\alpha)=\sqrt{1-\alpha^2 \|g_k\|^2}x_k+\alpha g_k, \ \ -\frac{1}{\|g_k\|}\le \alpha \le \frac{1}{\|g_k\|}.$$

If $\alpha_k \equiv \frac{1}{\|g_k\|}$, then SSPM method will reduce to the power method.
 For simplicity, if $\alpha_k \neq \frac{1}{\|g_k\|}$, the iterate can be expressed as $x(\alpha)=x_k+ \sigma g_k$ with $\sigma=\frac{\alpha}{\sqrt{1-\alpha^2 \|g_k\|^2}}$. In order to solve (\ref{2-dimenional subproblem}), one need to solve a equation like $\nabla f(x_k+\sigma g_k)^T g_k=(\mathcal{A}(x_k+\sigma g_k)^{m-1})^Tg_k= 0$. For each iteration, the computational cost of SSPM method is $m$ times than that of power method. As shown in \cite{HCD15}, the main computational cost of SSPM is the tensor-vector multiplications $\mathcal{A}x_k^{m-1}$ and $\mathcal{A}_k$ (defined in \cite{HCD15}), which requires $O(mn^m)$ operations and $O(m^2n^m)$ operations, respectively.

\section{Inexact gradient method}

%\noindent\rule{\textwidth}{1pt}
%\underline{\textbf{Algorithm 3: Gradient Method for (\ref{max-optimization-problem}) with optimal stepsize}} \\
%Given tensors $\mathcal{A} \in \mathbb{S}^{[m,n]}$ and $\mathcal{B}\in \mathbb{S}_+^{[m,n]}$, an initial unit iterate $x_0$. Let $\epsilon>0$ be the tolerance. Set k=0; Calculate gradient $g(x_k)$ and a search direction $d(x_k)$.\\% Let $\beta=1$ if we want to find the local maxima; otherwise, let $\beta=-1$.  \\
%\textbf{While} $\|g(x_k)\|> \epsilon$ \textbf{do}\\
%%calculate a search direction $d(x_k)$.\\
%\textbf{1:} Finding a stepsize $\frac{1}{\|g_k\|}\ge \alpha_k>0$ such that  $x(\alpha)=\sqrt{1-\alpha^2 \|g_k\|^2}x_k+\alpha g_k$ satisfying
%\begin{equation}\label{1-dimenional subproblem}
%\max f(x(\alpha))=\frac{\mathcal{A}x(\alpha)^m}{\mathcal{B}x(\alpha)^m}.
%\end{equation}
%\textbf{2:} Update the iterate $x_{k+1}=\frac{x(\alpha_k)}{\|x(\alpha_k)\|}$, calculate $g(x_{k+1})$.\\
%%\textbf{Step 3:} If $\|g(x_{k+1})\| \leq\epsilon $, then stop. Otherwise, set $k :=
%%k + 1$ and go to {Step 1}.\\
%\textbf{End while} \\
%\noindent\rule{\textwidth}{1pt}
Indicated by the idea in \cite{HCD15}, we can present a gradient method with optimal stepsize for computing the generalized tensor eigenpairs problem (\ref{max-optimization-problem}). But we don't want to present it here, since the computational cost per iterate is more expensive than Power method.
In this section, we firstly present the following inexact gradient method, and then establish its global convergence and linear convergence results under some suitable assumptions.

\noindent\rule{\textwidth}{1pt}
\underline{\textbf{Algorithm 3: Adaptive Gradient (AG) method for (\ref{max-optimization-problem})}} \\
Given tensors $\mathcal{A} \in \mathbb{S}^{[m,n]}$ and $\mathcal{B}\in \mathbb{S}_+^{[m,n]}$, an initial unit iterate $x_0$, parameter $\rho \in (0,1)$. Let $\epsilon>0$ be the tolerance. Set k=0; Calculate gradient $g(x_k)$ .\\% Let $\beta=1$ if we want to find the local maxima; otherwise, let $\beta=-1$.  \\
\textbf{While} $\|g(x_k)\|> \epsilon$ \textbf{do}\\
%calculate a search direction $d(x_k)$.\\
\textbf{1:} Generate a stepsize $\frac{1}{\|g_k\|}\ge \alpha_k>0$ such that  $x_k(\alpha)=\sqrt{1-\alpha^2 \|g_k\|^2}x_k+\alpha g_k$ satisfying
\begin{equation}\label{linesearch}
f(x_k(\alpha))\ge f(x_k)+\rho \alpha \|g(x_k)\|^2
\end{equation}
\textbf{2:} Update the iterate $x_{k+1}=x_k(\alpha_k)$, calculate $g(x_{k+1})$.\\
%\textbf{Step 3:} If $\|g(x_{k+1})\| \leq\epsilon $, then stop. Otherwise, set $k :=
%k + 1$ and go to {Step 1}.\\
\textbf{End while} \\
\noindent\rule{\textwidth}{1pt}

It is clear that $x_{k+1}\in \mathbb{S}^{n-1}\cap span\{x_k, g_k\}$. Moreover, by using (\ref{xtgx}), we can show the first-order gain per iterate is $g_k^T(x_{k+1}-x_k)=\alpha_k \|g_k\|^2$.  Since the spherical feasible region $\mathbb{S}^{n-1}$ is compact, $\mathcal{B}x^m$ is positive and bounds away from zero, we can get that all the functions and gradients of the objective (\ref{max-optimization-problem}) at feasible points are bounded \cite{ChenQi15}, i.e., there exists a constant $M>0$ such that for all $x\in \mathbb{S}^{n-1}$,
 \begin{equation}\label{boundofg}
 |f(x)|\le M, \ \ \mbox{and} \ \ \|g(x)\|\le M. %, \ \ \mbox{and} \ \ \|H(x)\|\le M.
 \end{equation}

 The following theorem indicates that the Algorithm 3 is convergent to the KKT point of the problem (\ref{max-optimization-problem}). The constructive proof is motivated by the idea in \cite{HCD15}.

\begin{theorem} \label{the:globalconvergence} Suppose that the gradient $g(x)$ is Lipschitz continuous on the unit shpere. Let $\{x_{k}\}$ is generated by Algorithm 3. Then the inexact curvilinear search condition defined in (\ref{linesearch}) is well-defined and there exists a positive constant $c>0$ such that
\begin{equation}\label{eq:boundofstepsize}
 f(x_{k+1})-f(x_k) \ge c \|g_k\|^2.
 \end{equation}
 Furthermore, $$\lim_{x\rightarrow\infty}\|g_k\|=0.$$
\end{theorem}

\noindent{\bf Proof.} Firstly, we have
\begin{equation}\label{eq:dxk}
x_k'(\alpha)=\frac{-\alpha\|g_k\|^2}{\sqrt{1-\alpha^2\|g_k\|^2}}x_k+g_k.
\end{equation}
Furthermore, we can obtain
\begin{equation}\label{eq:derivefxa}
\frac{\mbox{d}f(x_k(\alpha))}{\mbox{d}\alpha} |_{\alpha=0}=g(x_k(\alpha))^Tx_k'(\alpha)|_{\alpha=0}=\|g_k\|^2.
\end{equation}
Let $d_k(\alpha)=x_k(\alpha)-x_k$, using (\ref{xtgx}), we can derived that for any constant $\rho\in (0,1)$, there exists a positive scalar $\bar{\alpha} \le \frac{1}{\|g_k\|}$ such that for all
$\alpha \in (0, \bar{\alpha}]$,
$$f(x_k(\alpha))- f(x_k) \ge \rho \langle g_k, d_k(\alpha) \rangle = \rho \alpha \|g(x_k)\|^2.$$

Considering the gap between $f(x_k(\alpha))$ and $f(x_k)$, similarly to the proof of Lemma 4.2 in \cite{HCD15}, we can get
\begin{equation*}%\label{eq:gapfk}
\begin{split}
f(x_k(\alpha))-f(x_k)=&\int_0^{\alpha}g(x_k(t))^Tx'_k(t)dt\\
=&\alpha g(x_k(0))^Tx'_k(0)+\int_0^{\alpha}g(x_k(t))^T[x'_k(t)-x'_k(0)]dt\\
&+\int_0^{\alpha}[g(x_k(t))-g(x_k(0))]^Tx'_k(0)dt\\
&(\mbox{Using (\ref{eq:dxk}), (\ref{eq:derivefxa}), and Lipschitz condition})\\
\ge&\alpha \|g_k\|^2-M\int_0^{\alpha}\|x'_k(t)-x'_k(0)\|dt-L\|g_k\|\int_0^{\alpha}\|x_k(t)-x_k(0)\|dt\\
\ge&\alpha \|g_k\|^2-M\|g_k\|^2\int_0^{\alpha} \frac{t}{\sqrt{1-t^2\|g_k\|^2}}dt-L\|g_k\|\int_0^{\alpha}\sqrt{2-2\sqrt{1-t^2\|g_k\|^2}}dt\\
%\mbox{(\ref{ineq:A3 finite increase})}\leq&{\|{\bm u}^n-{\bm u}^\ast\|}^2_{{\bm M}_n}-{\|{\bm u}^{n+1}-{\bm u}^\ast\|}^2_{{\bm M}_{n+1}}\\
\end{split}
\end{equation*}
Without loss of generality, assume that $\alpha_k \le \tilde{\alpha}\triangleq\frac{\sqrt{3}}{2\|g_k\|}$, then for $t\le \alpha \le \tilde{\alpha}$, we have
$$\frac{1}{\sqrt{1-t^2\|g_k\|^2}}\le\frac{1}{\sqrt{1-\tilde{\alpha}^2\|g_k\|^2}}=2,$$
and
$$\sqrt{2-2\sqrt{1-t^2\|g_k\|^2}}\le 2t\|g_k\|.$$
So, we can obtain that
\begin{equation}\label{eq:gapfk1}
f(x_k(\alpha))-f(x_k)\ge (\alpha-M \alpha^2-L \alpha^2) \|g_k\|^2=(1-M\alpha-L\alpha)\alpha\|g_k\|^2.
\end{equation}

Set $\bar{\alpha}=\frac{1}{2(M+L)}$, we have $1>1-M\alpha-L\alpha\ge \frac{1}{2}$ for all $\alpha \in (0,\bar{\alpha}]$. It follows from (\ref{eq:gapfk1}) that (\ref{linesearch}) holds for all $\alpha \in (0,\bar{\alpha}]$ with $\rho=\frac{1}{2}$.
So, by using a backward strategy in curvilinear search, one can claim that the stepsize $\alpha_k$ is bounded from below. That is to say,
there exists a positive constant $c$ (e.g. $c=\frac{1}{4(M+L)}$) such that (\ref{eq:boundofstepsize}) holds. As $f(x)$ is bounded on the unit sphere, by (\ref{eq:boundofstepsize}), it is easy to prove that $\sum_{k}\|g_k\|^2< +\infty$. Therefore, $\lim_{x\rightarrow\infty}\|g_k\|=0;$ namely, the Algorithm 3 is globally convergent.    $\quad\Box$\\

In the rest of this section, we would like to establish the linear convergence rate of the Algorithm 3 under the assumption of second order sufficient condition. For convenience, rewrite (\ref{max-optimization-problem}) as
\begin{equation}\label{reformula-max-optimization-problem}
\max f(x)=\frac{\mathcal{A}x^m}{\mathcal{B}x^m} \ \ \mbox{subject to} \ \ c(x)=x^Tx-1=0.
\end{equation}
The Lagrangian function is $L(x,\mu)=f(x)-\mu c(x),$ and its gradient and Hessian are
$$\nabla_x L(x,\mu)=\nabla f(x)-\mu x, \mbox{and} \ \nabla_{x}^2 L(x,\mu)=\nabla^2 f(x)-\mu I_n.$$
At the KKT point $(x^*,\mu^*)$, we have $\nabla_x L(x^*,\mu^*)=0.$ So, $\nabla f(x^*)^Tx^*-\mu^* \|x^*\|^2=0$. By (\ref{xtgx}), we know $\mu^*=\nabla f(x^*)^Tx^*=0.$ We can formulate the sufficient conditions for $x^*$ being a strict local maximizer of (\ref{reformula-max-optimization-problem}) as: $g(x^*)=0$ and
\begin{equation}\label{eq:2ordercondition}
v^TH(x^*)v<0, \forall v \in \nabla c(x^*)^{\bot}\cap \mathbb{S}^{n-1},
\end{equation}
where $g(x)=\nabla f(x), \ H(x)=\nabla^2 f(x)$, and $x^{\bot}\equiv \{y\in \Re^n | x\perp y\}.$
By Theorem \ref{the:globalconvergence}, the sequence $\{x_k\}$ generated by Algorithm 3 is convergent to a KKT point $x^*$ with $g(x^*)=0$. If we further assume that the assumption of second sufficient condition (\ref{eq:2ordercondition}) holds at $x^*$, then the following linear convergence theorem could be established for Algorithm 3.

\begin{theorem}
 Let $\{x_{k}\}$ is generated by Algorithm 3. Suppose that the gradient $g(x)$ is Lipschitz continuous on the unit sphere and the second sufficient condition (\ref{eq:2ordercondition}) holds at the KKT point $x^*$. Then $\{f(x_k)\}$ converges to $f(x^*)$ linearly.
\end{theorem}
{\bf Proof.} In order to show $\{f(x_k)\}$ converges to $f(x^*)$ linearly, we need to prove
\begin{equation}
0<\lim_{k\rightarrow\infty}\frac{f(x^*)-f(x_{k+1})}{f(x^*)-f(x_{k})}=1-\lim_{k\rightarrow\infty}\frac{f(x_{k+1})-f(x_{k})}{f(x^*)-f(x_{k})}<1.
\end{equation}
To end of this, we firstly deduce that $f(x^*)-f(x_k)>0$. Project $x_k-x^*$ on the orthogonal space of $\nabla c(x^*)=2x^*$, we get
$$v_k=x_k-x^*-((x_k-x^*)^T\nabla c(x^*))\nabla c(x^*)=x_k-x^*-4((x_k-x^*)^Tx^*)x^*.$$
It is clear that $v_k \in \mathbb{S}^{n-1}$, since $\|x_k\|^2=\|x^*\|^2=1$. Notice that $-(x_k-x^*)^Tx^*=1-x_k^Tx^*=\frac{1}{2}(x_k-x^*)^T(x_k-x^*)$,
we can obtain that
$$x_k-x^*=v_k-2\|x_k-x^*\|^2x^*=v_k+O(\|x_k-x^*\|^2).$$

By the Taylor expansion, and the second order sufficient condition (\ref{eq:2ordercondition}), we have
\begin{equation*}%\label{eq:gapfk}
\begin{split}
f(x^*)-f(x_k)=&-\frac{1}{2}(x_k-x^*)^TH(x^*)(x_k-x^*)+o(\|x_k-x^*\|^2)\\
=&-\frac{1}{2}(v_k+O(\|x_k-x^*\|^2))^TH(x^*)(v_k+O(\|x_k-x^*\|^2))+o(\|x_k-x^*\|^2)\\
=&-\frac{1}{2}v_k^TH(x^*)v_k+o(\|x_k-x^*\|^2)>0.\\
%>&0.
\end{split}
\end{equation*}

Secondly, let us consider the following limitation:

\begin{equation*}%\label{eq:gapfk}
\begin{split}
\lim_{k\rightarrow\infty}\frac{\|g_k\|^2}{f(x^*)-f(x_k)}=&\lim_{k\rightarrow\infty}\frac{\|H(x^*)(x_k-x^*)\|^2}{-\frac{1}{2}v_k^TH(x^*)v_k}=\lim_{k\rightarrow\infty}\frac{2\|H(x^*)v_k\|^2}{-v_k^TH(x^*)v_k}.\\
%=&\lim_{k\rightarrow\infty}\frac{2\|H(x^*)v_k\|^2}{-v_k^TH(x^*)v_k}\\
%=&-\frac{1}{2}v_k^TH(x^*)v_k+o(\|x_k-x^*\|^2)\\
%>&0.
\end{split}
\end{equation*}
Let $U\in \Re^{n\times(n-1)}$ be the orthogonal complement of the vector $\nabla c(x^*)$. For any $v\in \nabla c(x^*)^{\bot}$, there exist $w\in \Re^{n-1}$ such that $v=Uw.$  Suppose that the Cholesky decomposition of the positive definite matrix $-U^TH(x^*)U$ is $D^TD$. Denoting $y=Dw$. Notice that $\|U\|_2=\sqrt{\|U^TU\|}=1,$ we can derive that:
\begin{equation*}%\label{eq:gapfk}
\begin{split}
\frac{2\|H(x^*)v\|^2}{-v^TH(x^*)v}=&\frac{2\|H(x^*)Uw\|^2}{-w^TU^TH(x^*)Uw}\\
=&\frac{2\|U\|\|H(x^*)Uw\|^2}{-w^TU^TH(x^*)Uw}\\
\ge&\frac{2\|U^TH(x^*)Uw\|^2}{-w^TU^TH(x^*)Uw}\\
=&\frac{-2y^TU^TH(x^*)Uy}{y^Ty}\\
\ge&2\lambda_{n-1},
\end{split}
\end{equation*}
where $\lambda_{n-1}$ is the smallest eigenvalue of the matrix $-U^TH(x^*)U$. Therefore, we have
\begin{equation}\label{eq:limofdeltaf}
\lim_{k\rightarrow\infty}\frac{\|g_k\|^2}{f(x^*)-f(x_k)}\ge 2\lambda_{n-1}.
\end{equation}
It follows from (\ref{eq:boundofstepsize}) and (\ref{eq:limofdeltaf}) that

\begin{equation*}%\label{eq:gapfk}
\begin{split}
\lim_{k\rightarrow\infty}\frac{f(x^*)-f(x_{k+1})}{f(x^*)-f(x_{k})}=&1-\lim_{k\rightarrow\infty}\frac{f(x_{k+1})-f(x_{k})}{f(x^*)-f(x_{k})}\\
\le&1-\lim_{k\rightarrow\infty}\frac{c\|g_k\|^2}{f(x^*)-f(x_{k})}\\
\le&1-2c\lambda_{n-1}\\
\le&1.
\end{split}
\end{equation*}
The proof is completed. $\quad\Box$\\

%%%%%%%%%%%%%%%%%%%%%%%%%%%%%%%%%%%%%%%%%
%%%%%%%%%%%%%%%%%%%%%%%%%%%%%%%%%%%%%%%%%%%%

\section{Numerical experiments}

In this section, we present some numerical results to illustrate the effectiveness of the proposed adaptive gradient (AG) method, which was compared
with the GEAP method--an adaptive shifted power method proposed by Tamara G. Kolda and Jackson R. Mayo \cite{KoldaMayo14}.
The experiments were done on a laptop with Intel Core 2 Duo CPU with a 4GB RAM, using MATLAB R2014b, and the Tensor Toolbox \cite{TensorToolbox}.
We set the parameter $\rho=0.001$ and initial guess of the stepsize $\alpha$ in (\ref{linesearch}) is generated by $\min (\frac{1}{\|g_k\|}, \frac{\|\triangle x_k\|}{\|\triangle g_k\|})$. If this initial guess can not satisfy the line search condition (\ref{linesearch}), then we truncate it as $\alpha=0.5*\alpha$, and try it again. Generally, once or twice is enough in our experiments.

In all numerical experiments, we stop the iterates once $|\lambda_{k+1}-\lambda_k|\le 10^{-10}$. The maximum iterations is 500.
\subsection{Comparison with GEAP for computing Z-eigenpairs}
The following example is originally from \cite{KR02} and was used in evaluating the SS-HOPM algorithm in \cite{KoldaMayo11} and the GEAP algorithm in \cite{KoldaMayo14} for computing Z-eigenpairs.

\begin{figure}\label{fig:AG-GEAP-Z-eigen}
$$
\begin{array}{cc}
\centerline{\includegraphics[width=1.1\textwidth]{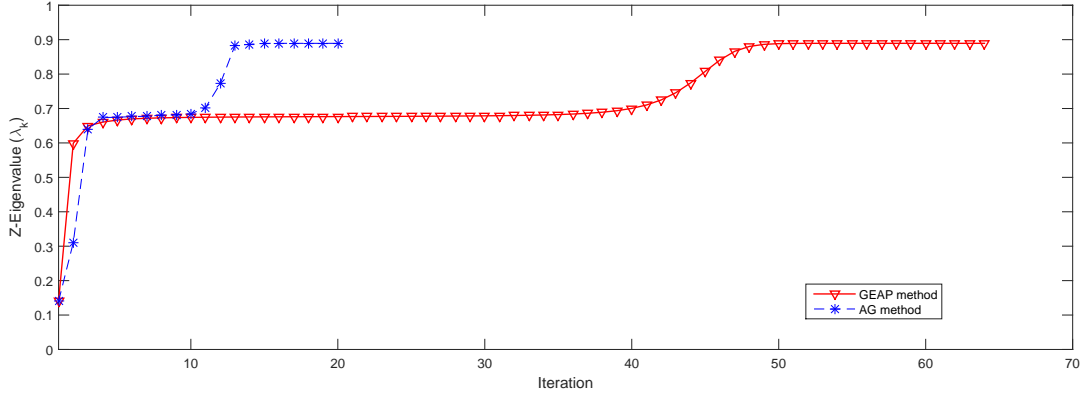}}
\end{array}
$$
 \caption{Comparison with GEAP algorithm for computing Z-eigenvalues of $\mathcal{A}$ from {\bf \it Example 1}, and the starting point is $x_0=[0.0417\  -0.5618\ \  0.6848]$ }

\end{figure}

{\it Example 1} (Kofidis and Regalia \cite{KR02}). Let $\mathcal{A} \in \mathbb{S}^{[4,3]}$ be the symmetric tensor defined by
\begin{equation*}%\label{eq:gapfk}
\begin{split}
a_{1111}=0.2883, \ \ a_{1112}=-0.0031,\ \  a_{1113}=0.1973, \ \ a_{1122}=-0.2485,&\\
a_{1223}= 0.1862, \ \ \  a_{1133}=\ 0.3847,\ \  a_{1222}=0.2972, \ \ a_{1123}=-0.2939,&\\
a_{1233}=0.0919, \ \ a_{1333}=-0.3619,\ \  a_{2222}=0.1241, \ \ a_{2223}= -0.3420,&\\
a_{2233}=0.2127, \ \ a_{2333}=0.2727,\ \  a_{3333}=-0.3054.&
\end{split}
\end{equation*}

To compare the convergence in terms of the number of iterations. Figure 1 shows the results for computing Z-eigenvalues of $\mathcal{A}$ from {\bf \it Example 1}, and the starting point is $x_0=[0.0417\  -0.5618\ \  0.6848]$. In this case, both of Adaptive Gradient (AG) method and GEAP method can find the largest Z-eigenvalue 0.8893. AG method just need run 19 iterations in 0.168521 seconds while GEAP method need run 63 iterations in  0.469648 seconds.

We used 1000 random starting guesses, each entry selected uniformly at random from the interval $[-1,1]$. For each set of experiments, the same set of random starts was used. For the largest eigenpair, we list the number of occurrences in the 1000 experiments. We also list the median number of iterations until convergence, the average error and the average run time in the 1000 experiments in Tables 1-4.
%The error is computed as $\|\mathcal{A}x^{m-1}-\lambda x\|_2$.
As we can see from Tables 1-4, Adaptive Gradient (AG) method is much faster than GEAP method and could reach the largest eigenpair with a higher probability.

\begin{center}
\small{Table 1. Comparison results for computing Z-eigenvalues of $\mathcal{A}$ from {\bf \it Example 1}.}\\
\begin{tabular}{|c|c|c|c|c|c|}
\hline
Alg. & Occ. &  $\lambda$ & Its.& Error& Time (sec.)\\
\hline
GEAP& 49.9\% &  0.8893 & 27.06 &  5.69e-11 & 0.1632 \\
AG& 56.6\% &  0.8893 & 13.81 &  1.74e-11 & 0.1205 \\
\hline
\end{tabular}
\end{center}

{\it Example 2} (Nie and Wang \cite{CDN14}). Let $\mathcal{A} \in \mathbb{S}^{[4,n]}$ be the symmetric tensor defined by
$$a_{ijkl}=\sin(i+j+k+l) \ \ (1\le i,j,k,l\le n).$$
For the case of $n=5$, there are five real Z-eigenvalues  which are respectively
$$\lambda_1=7.2595,\ \lambda_2=4.6408,\ \lambda_3=0.0000,\ \lambda_4=-3.9204,\ \lambda_5=-8.8463.$$

\begin{center}
\small{Table 2. Comparison results for computing Z-eigenvalues of $\mathcal{A}$ from {\bf \it Example 2}.}\\
\begin{tabular}{|c|c|c|c|c|c|}
\hline
Alg. & Occ. &  $\lambda$ & Its.& Error& Time (sec.)\\
\hline
GEAP& 48.2\% & 7.2595 & 50.01 &  7.68e-11 & 0.3235 \\
AG& 54.6\% &  7.2595 & 24.85 &  4.72e-11 & 0.2286 \\
\hline
\end{tabular}
\end{center}

{\it Example 3} (Nie and Wang \cite{CDN14}). Let $\mathcal{A} \in \mathbb{S}^{[4,n]}$ be the symmetric tensor defined by
$$a_{ijkl}=\tan(i)+\tan(j)+\tan(k)+\tan(l) \ \ (1\le i,j,k,l\le n).$$
For the case of $n=5$, there are five real Z-eigenvalues  which are respectively
$$\lambda_1=34.5304,\ \lambda_2=0.0000,\ \lambda_3=-101.1994.$$

\begin{center}
\small{Table 3. Comparison results for computing Z-eigenvalues of $\mathcal{A}$ from {\bf \it Example 3}.}\\
\begin{tabular}{|c|c|c|c|c|c|}
\hline
Alg. & Occ. &  $\lambda$ & Its.& Error& Time (sec.)\\
\hline
GEAP& 64.0\% &  34.5304 & 28.07 &  5.84e-11 & 0.1701 \\
AG& 83.9\% &  34.5304 & 17.70 &  3.17e-11 & 0.1544 \\
\hline
\end{tabular}
\end{center}

{\it Example 4} (Nie and Wang \cite{CDN14}). Let $\mathcal{A} \in \mathbb{S}^{[4,n]}$ be the symmetric tensor defined by
$$a_{ijkl}=\arctan((-1)^{i}\frac{i}{n})+\arctan((-1)^{j}\frac{j}{n})+\arctan((-1)^{k}\frac{k}{n})+\arctan((-1)^{l}\frac{l}{n}) $$
%For the case of $n=5$, there are five real Z-eigenvalues  which are respectively
%$$\lambda_1=34.5304,\ \lambda_2=0.0000,\ \lambda_3=-101.1994.$$

\begin{center}
\small{Table 4. Comparison results for computing Z-eigenvalues of $\mathcal{A}$ from {\bf \it Example 4} (n=5).}\\
\begin{tabular}{|c|c|c|c|c|c|}
\hline
Alg. & Occ. &  $\lambda$ & Its.& Error& Time (sec.)\\
\hline
GEAP& 65.2\% & 13.0779 & 22.43 &  5.64e-11 & 0.1453 \\
AG& 87.7\% &  13.0779 & 13.88 &  2.66e-11 & 0.1242 \\
\hline
\end{tabular}
\end{center}

\subsection{Comparison with GEAP for computing H-eigenpairs}

In this subsection, we test the proposed AG method with comparison to GEAP method on finding H-eigenpairs.

To compare the convergence in terms of the number of iterations. Figure 2 shows the results for computing H-eigenvalues of $\mathcal{A}$ from {\bf \it Example 5}(n=5), and the starting point is $x_0=[-0.8181 \; -0.4264 \;-0.0163 \;\;0.1198 \;-0.1574]$. In this case, GEAP fails to stop in 500 iterations. But Adaptive Gradient (AG) method can find the largest H-eigenvalue 0.8 after running 25 iterations in 0.2183 seconds.

{\it Example 5}.  Let $\mathcal{A} \in \mathbb{S}^{[4,n]}$ be the diagonal tensor defined by
$a_{iiii}=\frac{i-1}{i}. $

\begin{center}
\small{Table 5. Comparison results for computing H-eigenvalues of $\mathcal{A}$ from {\bf \it Example 5} (n=5).}\\
\begin{tabular}{|c|c|c|c|c|c|}
\hline
Alg. & Occ. &  $\lambda$ & Its.& Error& Time (sec.)\\
\hline
GEAP& 81\% & 0.8 & 500 &  4.04e-05 & 2.6941 \\
AG& 94\% &  0.8 & 14.48 &  4.34e-11 & 0.0872 \\
\hline
\end{tabular}
\end{center}

\begin{figure}\label{fig:AG-GEAP-Z-eigen}
$$
\begin{array}{cc}
\centerline{\includegraphics[width=1.1\textwidth]{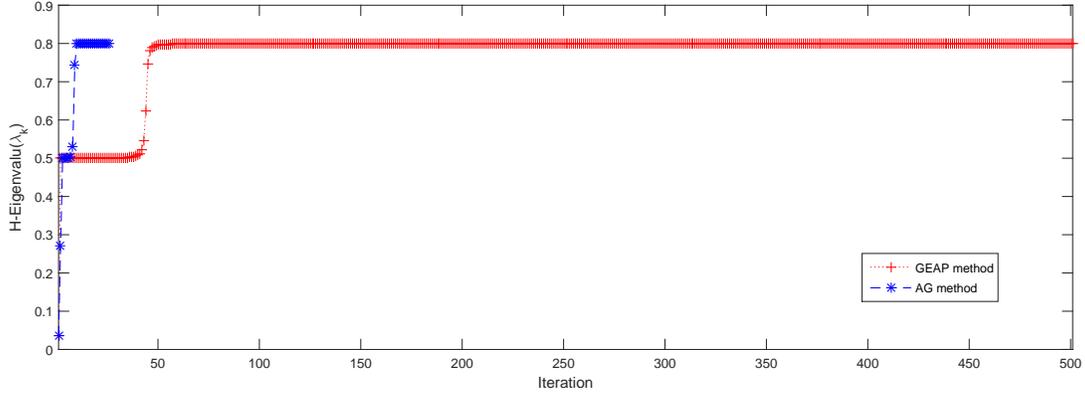}}
\end{array}
$$
 \caption{Comparison with GEAP algorithm for computing H-eigenvalues of $\mathcal{A}$ from {\bf \it Example 5} (n=5), and the starting point is $x_0=[-0.8181 \; -0.4264 \;-0.0163 \;0.1198 \;-0.1574]$ }

\end{figure}

\begin{figure}\label{fig:AG-GEAP-Z-eigen}
$$
\begin{array}{cc}
\centerline{\includegraphics[width=1.1\textwidth]{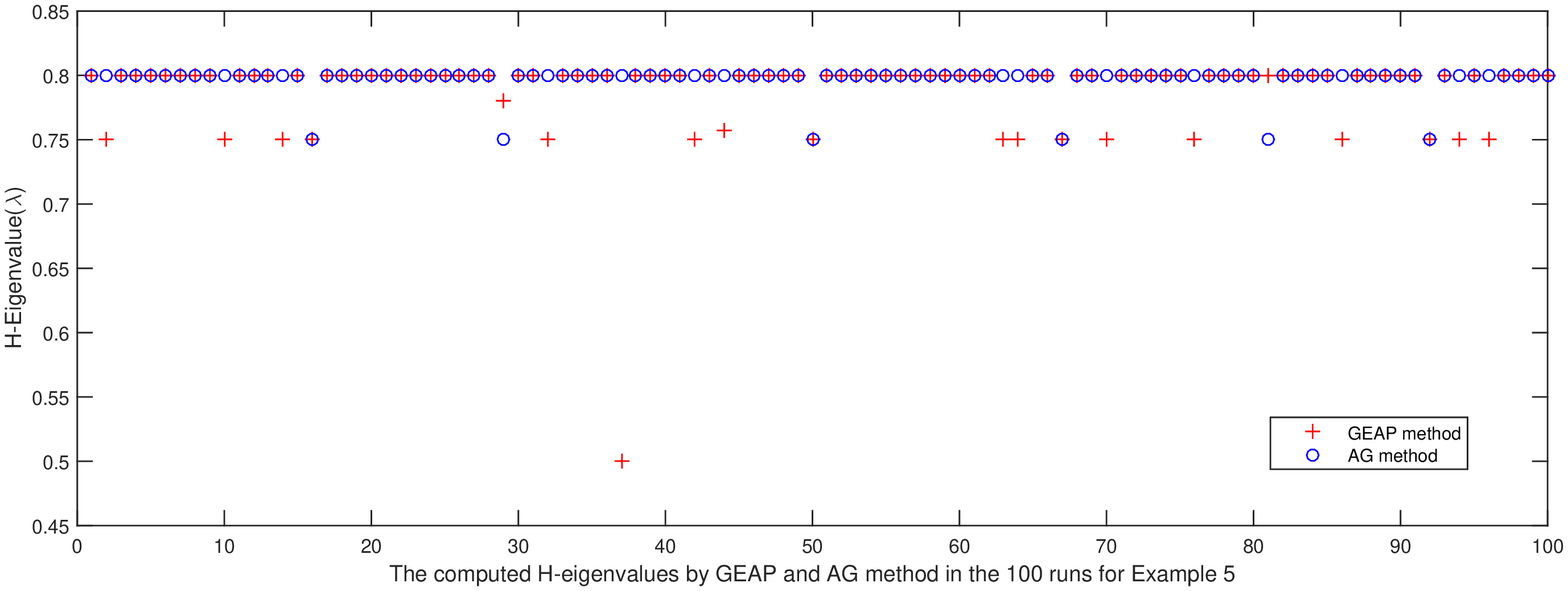}}
\end{array}
$$
 \caption{The computed H-eigenvalues by GEAP and AG method in the 100 runs on the $\mathcal{A}$ from {\bf \it Example 5} (n=5). }

\end{figure}

{\it Example 6} (Nie and Wang \cite{CDN14}).  Let $\mathcal{A} \in \mathbb{S}^{[4,n]}$ be the tensor defined by
$$a_{ijkl}=\frac{(-1)^i}{i}+\frac{(-1)^j}{j}+\frac{(-1)^k}{k}+\frac{(-1)^l}{l}, (1\le i,j,k,l \le n). $$

\begin{center}
\small{Table 6. Comparison results for computing H-eigenvalues of $\mathcal{A}$ from {\bf \it Example 6} (n=5).}\\
\begin{tabular}{|c|c|c|c|c|c|}
\hline
Alg. & Occ. &  $\lambda$ & Its.& Error& Time (sec.)\\
\hline
GEAP& 61\% & 34.3676 & 20.94 &  5.32e-11 & 0.1209 \\
AG& 100\% &  34.3676 & 15.71 &  1.93e-11 & 0.2650 \\
\hline
\end{tabular}
\end{center}

\begin{figure}\label{fig:AG-GEAP-Z-eigen}
$$
\begin{array}{cc}
\centerline{\includegraphics[width=1.1\textwidth]{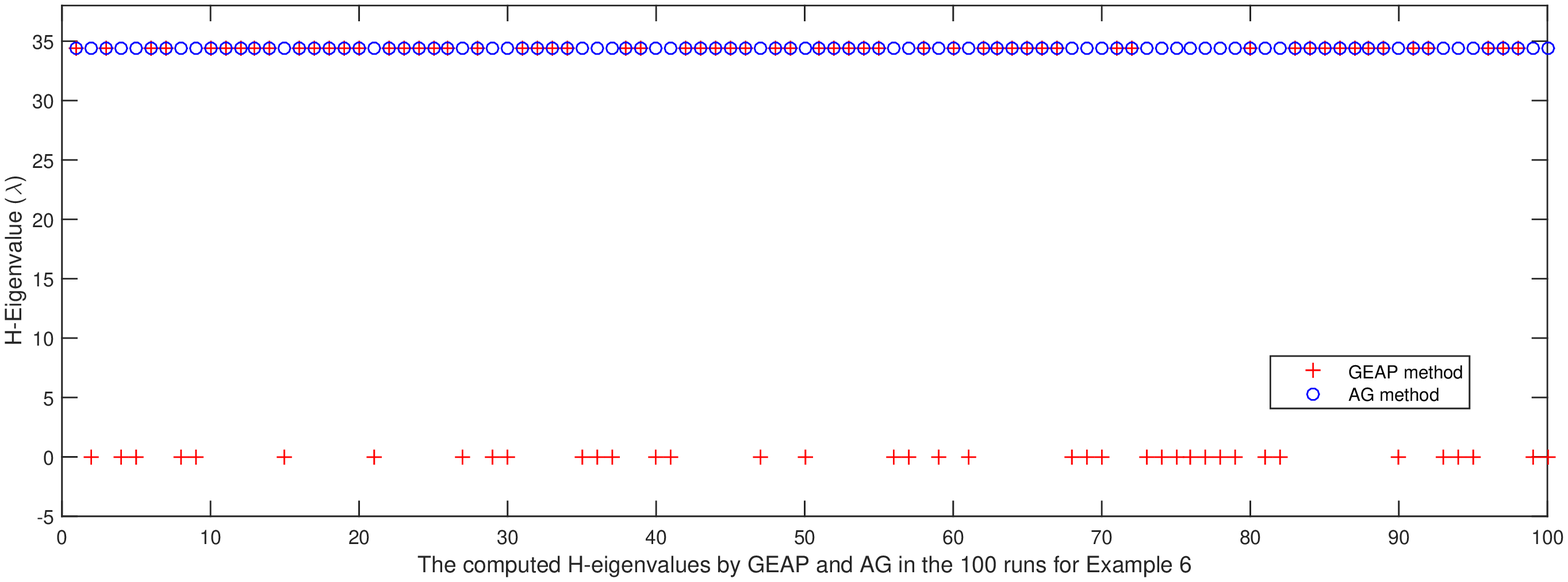}}
\end{array}
$$
 \caption{The computed H-eigenvalues by GEAP and AG method in the 100 runs on the $\mathcal{A}$ from {\bf \it Example 6} (n=5). }

\end{figure}

{\it Example 7}.  Let $\mathcal{A} \in \mathbb{S}^{[4,n]}$ be the tensor defined by
$$a_{iiii}=2i, (1\le i\le n) \ \mbox{and} \ \ a_{1123}=4b, $$
 here $b$ is a parameter. Then use the {\bf symmetrize} function in the Matlab Tensor Toolbox to symmetrize it.

\begin{center}
\small{Table 7. Comparison results for computing H-eigenvalues of $\mathcal{A}$ from {\bf \it Example 7} ($n=3,b=1$).}\\
\begin{tabular}{|c|c|c|c|c|c|}
\hline
Alg. & Occ. &  $\lambda$ & Its.& Error& Time (sec.)\\
\hline
GEAP& 71\% & 6.112 & 179.01 &  1.99e-08 & 0.9536 \\
AG& 100\% &  6.112 & 50.52 &  5.94e-11 & 0.5198 \\
\hline
\end{tabular}
\end{center}

\begin{figure}\label{fig:AG-GEAP-Z-eigen}
$$
\begin{array}{cc}
\centerline{\includegraphics[width=1.1\textwidth]{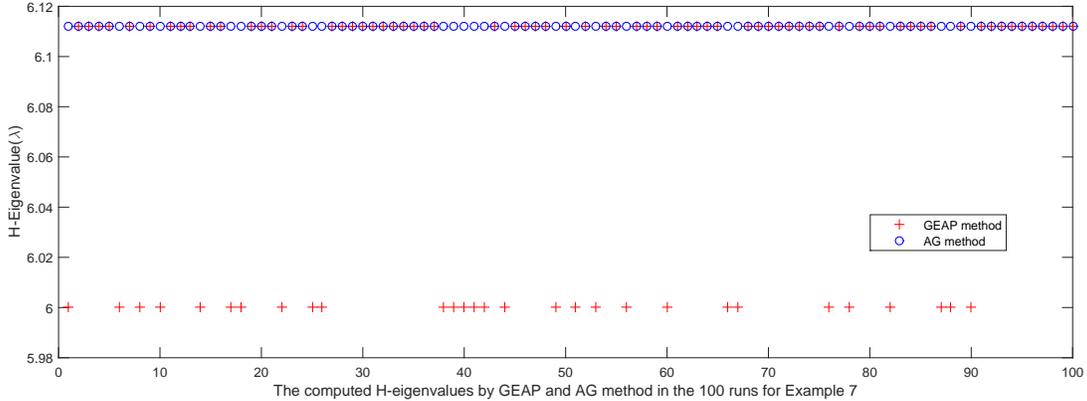}}
\end{array}
$$
 \caption{The computed H-eigenvalues by GEAP and AG method in the 100 runs on the $\mathcal{A}$ from {\bf \it Example 7} (n=3,b=1). }

\end{figure}

We used 100 random starting guesses to test AG method and GEAP method for computing H-eigenvalues of $\mathcal{A}$ from {\bf \it Examples 5-7}. For each set of experiments, the same set of random starts was used. For the largest eigenpair, we list the number of occurrences in the 100 experiments. We also list the median number of iterations until convergence, the average error and the average run time in the 100 experiments in Tables 5-7. As we can see from Table 5, GEAP method fails to stop in 500 iterations for all of the 100 test experiments for Example 5. But GEAP can slowly approach to the largest H-eigenvalue in 81 test experiments as shown in Figure 3.  Adaptive Gradient (AG) method much faster than GEAP method and could reach the largest eigenpair with a higher probability. Especially, for Examples 6 and 7, Adaptive Gradient (AG) method could find the largest H-eigenvalue in all of the 100 experiments.

%%%%%%%%%%%%%%%%%%%%%%%%%%%%%%%%%%%%%%%%%%%%

\section{Conclusion}

In this paper, we introduced an adaptive gradient (AG) method for generalized tensor eigenpairs, which could be viewed as an inexact version of the gradient method with optimal stepsize for finding Z-eigenvalues of tensor in \cite{HCD15}. What we have done is to use an inexact curvilinear search condition to replace the constraint on optimal stepsize. So, the computational complexity of AG method is much cheaper than SSPM method in \cite{HCD15}. Global convergence and linear convergence rate are established for the AG method for computing generalized eigenpairs of symmetric tensor. Some numerical experiments illustrated that the AG method is faster than GEAP and could reach the largest eigenpair with a higher probability.

\section*{Acknowledgements}

 This work was supported in part by the National Natural Science Foundation of China (No.61262026, 11571905, 11501100), NCET Programm of the Ministry of Education (NCET 13-0738), JGZX programm of Jiangxi Province (20112BCB23027), Natural Science Foundation of Jiangxi Province (20132BAB201026), science and technology programm of Jiangxi Education Committee (LDJH12088), Program for Innovative Research Team in University of Henan Province (14IRTSTHN023).


\begin{thebibliography}{99}




%\bibitem{BB88} J. Barzilai and J.M. Borwein, ``Two-point step size gradient methods", {\sl IMA J. Numer. Anal.}, 8 (1988) 141-148.



 \bibitem{TensorToolbox}   B. W. Bader, T. G. Kolda and others. MATLAB Tensor Toolbox Version 2.6, Available online, February 2015.  http://www.sandia.gov/~tgkolda/TensorToolbox/

 \bibitem{BV08} L. Bloy, R. Verma, ``On Computing the Underlying Fiber Directions from the Diffusion Orientation Distribution
Function", In Medical Image Computing and Computer-Assisted Intervention MICCAI, Vol. 2008. Springer: Berlin /
Heidelberg, 2008: 1-8.


 \bibitem{CetinUnal15} S. Cetin and G. Unal, ``A Higher-Order Tensor Vessel Tractography for Segmentation of Vascular Structures", {\sl IEEE Transactions on Medical Imaging} 34 (2015) 2172-2185.


 \bibitem{CPZ08}  K.C. Chang, K.J. Pearson and T. Zhang, ``Perron-Frobenius Theorem for nonnegative tensors",
 {\sl Communications in Mathematical Sciences}, 6 (2008) 507-520.

\bibitem{CPZ09}  K.C. Chang, K.J. Pearson and T. Zhang, ``On eigenvalue problems of real symmetric tensors",
 {\sl J. Math. Anal. Appl.}, 350 (2009) 416-422.

 \bibitem{CPZ11} K. Chang, K. Pearson, T. Zhang, ``Primitivity, the convergence of the NQZ method, and the largest eigenvalue for
nonnegative tensors",  {\sl SIAM Journal on Matrix Analysis and Applications} 32 (2011) 806-819.

\bibitem{CQZ13} K. Chang, L. Qi, and T. Zhang, ``A survey on the spectral theory of nonnegative tensors," {\sl Numerical Linear Algebra with Applications}, 20 (2013) 891-912.

\bibitem{CZ13} K. Chang, T. Zhang, ``On the uniqueness and nonuniqueness of the Z-eigenvector for transition probability tensors",
{\sl Journal of Mathematical Analysis and Applications} 408 (2013) 525-540.


\bibitem{CDHS13} Y. Chen, Y. Dai, D. Han, and W. Sun, ``Positive semidefinite generalized diffusion tensor
imaging via quadratic semidefinite programming", {\sl SIAM J. Imaging Sci.}, 6 (2013) 1531-1552.

\bibitem{ChenQi15}  Y. Chen, L. Qi and Q. Wang, ``Computing eigenvalues of large scale Hankel tensors",
 {\sl arXiv: 1504.07413v2}, May 2015.


\bibitem{CW14} M. Chu, S. Wu, ``On the second dominant eigenvalue affecting the Power method for transition probability tensors", Manuscript, 2014.

\bibitem{CDN14} C. Cui, Y. Dai and J. Nie, ``All real eigenvalues of symmetric tensors", {\sl SIAM J. Matrix Anal. Appl.}, 35 (2014) 1582-1601.

\bibitem{DW15} W. Ding and Y. Wei, ``Generalized tensor eigenvalue problems", {\sl SIAM J. Matrix Anal. Appl.}, 36 (2015) 1073-1099.





\bibitem{Han2013} L. Han, ``An unconstrained optimization approach for finding real eigenvalues of even order symmetric tensors'',{\sl Numer. Algebr. Control Optim.} 3 (2013) 583-599.



\bibitem{HCD15} C. Hao, C. Cui and Y. Dai, ``A sequential subspace projection method for extreme
Z-eigenvalues of supersymmetric tensors", {\sl Numer. Linear Algebr. Appl.}, 22 (2015) 283-298.

\bibitem{HCD15-PJO} C. Hao, C. Cui and Y. Dai, ``A feasible trust-region method for calculating extreme
Z-eigenvalues of symmetric tensors", {\sl Pacific J. Optim.}, 11 (2015) 291-307.


\bibitem{HillarLim2013} C. Hillar and L. Lim, ``Most tensor problems are NP-hard", {\sl J. ACM} 60 (2013) article No. 45: 1-39.


\bibitem{HLQ13} S. Hu, G. Li, L. Qi, Y. Song, ``Finding the maximum eigenvalue of essentially nonnegative symmetric tensors via sum of
squares programming on the largest eigenvalue of a symmetric nonnegative tensor", {\sl Journal of Optimization Theory
and Applications} 2013.



\bibitem{HHQ13} S. Hu, Z. Huang and L. Qi, ``Finding the extreme Z-eigenvalues of tensors via a sequential
SDPs method", {\sl Numer. Linear Algebr. Appl.}, 20 (2013) 972-984.


\bibitem{KR02} E. Kofidis and P. Regalia, ``On the best rank-1 approximation of higher-order supersymmetric
tensors", {\sl SIAM J. Matrix Anal. Appl.}, 23 (2002) 863-884.

\bibitem{KoldaMayo14} T. G. Kolda and J. R. Mayo, ``An adaptive shifted power methods for computing generalized tensor eigenpairs'', {\sl SIAM J. Matrix Anal. Appl.}, 35 (2014) 1563-1581.

\bibitem{KoldaMayo11} T. G. Kolda and J. R. Mayo, ``Shifted power method for computing tensor eigenpairs'', {\sl SIAM J. Matrix Anal. Appl.},  32 (2011) 1095-1124.




\bibitem{LMV00} L. De Lathauwer, B. De Moor, and J. Vandewalle, ``On the best rank-1 and rank-(R1,R2, . . . ,RN) approximation of higher-order tensors", {\sl SIAM J. Matrix Anal. Appl.}, 21 (2000)  1324-1342.

\bibitem{LQY13} G. Li, L. Qi, G. Yu, ``The Z-eigenvalues of a symmetric tensor and its application to spectral hypergraph theory",
{\sl Numerical Linear Algebra with Applications}, 20 (2013) 1001-1029.

\bibitem{LQY13-LAA}G. Li, L. Qi, G. Yu, ``Semismoothness of the maximum eigenvalue function of a symmetric tensor and its application", {\sl Linear Algebra and its Applications}, 438 (2013) 813-833.

\bibitem{LN13} W. Li, M. Ng, ``On the limiting probability distribution of a transition probability tensor", {\sl Linear and Multilinear Algebra}, 62 (2014) 362-385.

 \bibitem{Lim05} L.H. Lim, `` Singular values and eigenvalues of tensors, A variational approach",
 Proc. 1st IEEE International workshop on computational advances of multi-tensor adaptive processing, 2005, 129-132.

%\bibitem{LZI10} Y. Liu, G. Zhou, N. F. Ibrahim, ``An always convergent algorithm for the largest
%eigenvalue of an irreducible nonnegative tensor",  {\sl J. Comput Appl Math}, 235 (2010) 286-292.




\bibitem{NQZ10} M. Ng, L. Qi, G. Zhou, ``Finding the largest eigenvalue of a nonnegative tensor", {\sl SIAM Journal on Matrix Analysis and
Applications} 31 (2010) 1090-1099.

\bibitem{NQW08} Q. Ni, L. Qi and F. Wang, ``An eigenvalue method for testing positive definiteness of a
multivariate form", {\sl IEEE T. Automat. Contr.}, 53 (2008) 1096-1107.

\bibitem{NQ15-GlobalOpt} Q. Ni and L. Qi, ``A quadratically convergent algorithm for finding the largest eigenvalue of a nonnegative homogeneous polynomial map", {\sl Journal of Global Optimization}, 61 (2015)  627-641.


\bibitem{NQB14} G. Ni, L. Qi and M. Bai, ``Geometric measure of entanglement and U-eigenvalues of
tensors", {\sl SIAM J. Matrix Anal. Appl.}, 35 (2014) 73-87.


\bibitem{PLV05} J.M. Papy, L. De Lathauwer and S. Van Huffel, ``Exponential data fitting using multilinear
algebra: the single-channel and multi-channel case", {\sl Numer. Linear Algebr. Appl.}, 12 (2005) 809-826.

\bibitem{PLV09} J.M. Papy, L. De Lathauwer and S. Van Huffel, ``Exponential data fitting using multilinear
algebra: the decimative case", {\sl J. Chemometr.} 23 (2009) 341-351.


\bibitem{Qi05} L. Qi, ``Eigenvalues of a real supersymmetric tensor",
{\it J. Symbolic Computation,} 40 (2005) 1302-1324.

\bibitem{Qi07} L. Qi, ``Eigenvalues and invariants of tensor",
{\it J. Math. Anal. Appl.,} 325 (2007) 1363-1377.

\bibitem{QHW09} L. Qi, D. Han and E.X. Wu, ``Principal invariants
and inherent parameters of diffusion kurtosis tensors'', {\sl
Journal of Mathematical Analysis \& Applications}, 349 (2009)
165-180.

\bibitem{QWW} L. Qi, F. Wang and Y. Wang, ``Z-eigenvalue methods
for a global polynomial optimization problem'', {\sl Mathematical
Programming}, 118 (2009) 301-316.

\bibitem{QWW08} L. Qi, Y. Wang and E.X. Wu, ``D-Eigenvalues of diffusion kurtosis tensors'',
{\sl Journal of Computational and Applied Mathematics}, 221 (2008)
150-157.

\bibitem{QYW} L. Qi, G. Yu and E.X. Wu, ``Higher order positive semi-definite diffusion tensor
imaging'', {\sl SIAM Journal on Imaging Sciences}, 3 (2010) 416-433.

\bibitem{QYX} L. Qi, G. Yu and Y. Xu, ``Nonnegative diffusion orientation distribution function'', {\sl Journal of Mathematical Imaging and Vision}, 45 (2013) 103-113.

\bibitem{SS08} T. Schultz and H.-P. Seidel, ``Estimating crossing fibers: a tensor decomposition approach",
{\sl IEEE T. Vis. Comput. Gr.}, 14 (2008) 1635-1642.



\bibitem{SQ15} Y. Song and L. Qi, ``Tensor Complementarity Problem and Semi-positive Tensors", {\sl Journal of Optimization Theory and Applications}, to appear.

\bibitem{SY15} Y. Song and G. Yu, ``Properties of Solution Set of Tensor Complementarity Problem", {\sl arXiv:1508.00069v2}, Aug. 2015.

\bibitem{SJY08} L. Sun, S. Ji and J. Ye, ``Hypergraph spectral learning for multi-label classification", in Proceedings
of the 14th ACM SIGKDD International Conference on Knowledge Discovery and
Data Mining, ACM, 2008, pp. 668-676.


%
%\bibitem{ZCTW12} G. Zhou, L. Caccetta, K. Teo, S. Wu, ``Positive polynomial optimization over unit spheres and convex programming
%relaxations", {\sl SIAM Journal on Optimization} 22 (2012) 987-1008.

\bibitem{YYZ14} L. Yang, Q. Yang, X. Zhao, ``Quadratic third-order tensor optimization problem with quadratic constraints", {\sl Statistics, Optimization and Information Computing}, 2 (2014) 130-146.

\bibitem{YY11} Q. Yang, Y. Yang, ``Further results for Perron-Frobenius Theorem for nonnegative tensors II", {\sl SIAM Journal on Matrix
Analysis and Applications} 32 (2011) 1236-1250.

\bibitem{YY10} Y. Yang, Q. Yang, ``Further results for Perron-Frobenius Theorem for nonnegative tensors", {\sl SIAM Journal on Matrix
Analysis and Applications} 31 (2010) 2517-2530.

\bibitem{XC13}  J. Xie, A. Chang, ``H-eigenvalues of the signless Laplacian tensor for an even uniform hypergraph",  {\sl Frontiers of
Mathematics in China}, 8 (2013) 107-127.

\bibitem{ZN15-PJO} M. Zeng and Q. Ni, ``Quasi-Newton method for computing Z-eigenpairs of a symmetric tensor", {\sl Pacific J. Optim.},
11 (2015) 279-290.

\bibitem{ZZP13} F. Zhang, B. Zhou,  L. Peng, ``Gradient skewness tensors and local illumination detection for images", {\sl Journal of
Computational and Applied Mathematics} 237 (2013) 663-671.

\bibitem{ZQ12} L. Zhang, L. Qi, ``Linear convergence of an algorithm for computing the largest eigenvalue of a nonnegative tensor",
{\sl Numerical Linear Algebra with Applications} 19 (2012) 830-841.

\bibitem{ZQX12} L. Zhang, L. Qi, Y. Xu, ``Finding the largest eigenvalue of an irreducible nonnegative
tensor and linear convergence for weakly positive tensors",  {\sl J Comput Math}, 30 (2012) 24-33.




\end{thebibliography}
\end{document}